\documentclass{amsart}

\newtheorem{theorem}{Theorem}[section]

\theoremstyle{definition}

\theoremstyle{remark}
\newtheorem{remark}[theorem]{Remark}

\numberwithin{equation}{section}

\usepackage{amsmath}
\usepackage{amsfonts}
\usepackage{amssymb}
\usepackage{hyperref}
\usepackage[dvipsnames]{xcolor}
\usepackage{amsmath,systeme}

\begin{document}
\title[Rotational Rigidity]{A Rigidity Result for minimal Rotation hypersurfaces of 5d spaces of constant curvature}
\author{Aaron J. Tyrrell}
\date{April 2023}
\maketitle
\begin{abstract}
In this paper we show that a particular extrinsic pointwise hypersurface invariant is always non-positive on minimal hypersurfaces of constant curvature spaces and vanishes identically if and only if the hypersurface is rotational. We show this for hypersurfaces of 5-dimensional spaces of constant curvature but we conjecture that this should generalize to a similar result in other dimensions. 

\end{abstract}
\section{introduction}
Let $M^5(c)$ be a simply connected, complete Riemannian manifold with constant curvature $c,$ where $c$ is a real number.  
Do Carmo and Dajczer \cite{do2012rotation} studied rotational hypersurfaces of such manifolds and a consequence of Corollary $4.4$ of their work is that a non-totally geodesic minimal hypersurface of $M^5(c)$ is rotational if and only if it has exactly $3$ principal curvatures which are equivalent. Some other work on rotational hypersurfaces can be found here \cite{wang2018simons}, \cite{levitt1985symmetry}.
\section{Statement and proof of result}
\begin{theorem}[\textbf{Rotational Rigidity}]
Let $N$ be a minimal hypersurface of $M^5(c).$ Let $A$ be the second fundamental form of $N.$ Then define
\begin{equation}
S(N):=|A^2|^2-\frac{7}{12}|A|^4,
\end{equation}
then
\begin{equation*}
S(N)(p)\leq 0, \quad\forall\hspace{1mm} {p\in N}
\end{equation*}
with 
\begin{equation*}
S(N)\equiv 0
\end{equation*}
if and only if $N$ is a rotational hypersurface.
\end{theorem}
\begin{remark}
This gives another way to view minimal rotation hypersurfaces of $M^5(c).$ This is as minimizers amongst minimal hypersurfaces of the energy 
\begin{equation}
E_{rot}[N]:=\int_N \frac{7}{12}|A|^4-|A^2|^2dA,
\end{equation}
where $dA$ is the area form of the induced metric on $N.$
\end{remark}

\begin{proof}
Given $p\in N,$ let $\lambda_1,\lambda_2,\lambda_3,\lambda_4$ be the principal curvatures of $N$ at $p.$ Then letting
\begin{align}
SN(p)
&=\lambda_1^4+\lambda_2^4+\lambda_3^4+\lambda_4^4-\frac{7}{12}(\lambda_1^2+\lambda_2^2+\lambda_3^2+\lambda_4^2)^2\\\nonumber
&=:s(\lambda_1,\lambda_2,\lambda_3,\lambda_4).
\end{align}
We can apply Lagrange multipliers to $s(x,y,z,w)$ subject to the constraint $x+y+z+w=0$ to get the system:
\begin{align}
4x^3-\frac{28}{12}||x||^2x&=\lambda\\
4y£^3-\frac{28}{12}||x||^2y&=\lambda\\
4z^3-\frac{28}{12}||x||^2z&=\lambda\\
4w^3-\frac{28}{12}||x||^2w&=\lambda\\
x+y+z+w&=0,
\end{align}
where $||x||^2:=x^2+y^2+z^2+w^2.$ We also define $x_1:=x, x_2:=y, x_3:=z, x_4:=w.$\par
By the symmetry of these equations we see that the set of points satisfying these equations simultaneously is preserved under reflections through the hyperplanes $\{x_i=x_j : \hspace{1mm} i\neq j\}$. A natural approach is to set $x=y=z$ and try to solve.  
This yields $w=-3x, z=y,x=y$ and
\begin{equation*}-108x^3+84x^3=\lambda,
\end{equation*}
\begin{equation*}
x^3=-\frac{\lambda}{24},
\end{equation*}
and so
\begin{equation*}(x,y,z,w)=-\sqrt[3]{\lambda}\big(\frac{1}{\sqrt[3]{24}},\frac{1}{\sqrt[3]{24}},\frac{1}{\sqrt[3]{24}},-\frac{\sqrt[3]{9}}{\sqrt[3]{8}} \big).\end{equation*}
By the scale invariance of $(2.4)-(2.8)$ we can let $\lambda=24\alpha^3$ where $\alpha\in \mathbb{R}$ and we get critical points 
\begin{equation*}\alpha(-1,-1,-1,3) \hspace{2mm}\forall\hspace{1mm}\alpha\in\mathbb{R}.
\end{equation*}
Observe
\begin{equation*}s(\alpha(-1,-1,-1,3))=0.\end{equation*}
Clearly we could have chosen any one of $x,y,z,w$ to not be equal to the others, and therefore we get a critical point at $\alpha\sigma(-1,-1,-1,3)$ for any $\alpha \in\mathbb{R}$ and $\sigma\in S_4$ where $S_4$ is the symmetric group on 4 letters. 
Now to see that $0$ is the largest critical value, first consider that at a critical point $(x,y,z,w)$ subject to our constraint with $x,y,z$ distinct, $(2.4)$ and $(2.5)$ are equal so we must have
\begin{equation}
4x^3-\frac{28}{12}||x||^2x
=4y^3-\frac{28}{12}||x||^2y, 
\end{equation}
therefore
\begin{align}
4(x^3-y^3)&=\frac{28}{12}||x||^2(x-y). 
\end{align}
Factoring the left hand-side of $(2.10)$ gives 
\begin{align*}
4(x-y)(x^2+xy+y^2)&=\frac{28}{12}||x||^2(x-y), 
\end{align*}
recall $x\neq y$ so 
\begin{equation}
(x^2+xy+y^2)=\frac{7}{12}(x^2+y^2+z^2+w^2).
\end{equation}
The right hand-side of $(2.11)$ is independent of the choice of equations we chose to set equal in $(2.9)$ and since $x\neq z$ we get 
\begin{equation*}
x^2+xy+y^2=x^2+xz+z^2
\end{equation*}
this implies
\begin{equation*}
x(y-z)+(y-z)(y+z)=0,
\end{equation*}
and therefore
\begin{equation*}
(y-z)(x+y+z)=0
\end{equation*}
thus, since $y\neq z$ we get 
\begin{equation}
w=0.
\end{equation}
This implies that  $\lambda=0,$
now if any two of $x,y,z $ is non-zero then it follows that they differ by a sign. To see this take $(2.4)$ and see that
\begin{equation}
4x^2=\frac{28}{12}||x||^2,
\end{equation}
since the right hand-side of $(2.13)$ is independent of which of $(2.4)$, $(2.5),$ $(2.6)$ and $(2.7)$ we chose we get that $x_i^2=x_j^2$ if $x_i\neq 0$ and $x_j\neq 0.$
The only possibility is that two of $x, y, z$  are additive inverses and the other is $0.$ WLOG assume $x=-y$ then \begin{align}s(x,-x,0,0)&=2x^4-\frac{7}{3}x^4\\\nonumber
&=-\frac{1}{3}x^4.\end{align}Therefore any constrained critical point with 3 distinct coordinates is not a global maximum (we only considered the case where $x,y$ and $z$ are distinct but by symmetry this covers every case where three of the coordinates are distinct).
\par We are left to consider points of the form 
$\sigma(x,y,x,y)$ for $\sigma\in S_4$. Now by $(2.8)$ we must have $x=-y$ which implies the point has the form $(x,-x,x,-x).$
Clearly then 
\begin{align*}s(x,-x,x,-x)&=4x^4-\frac{28}{3}x^4\\
&=-\frac{26}{3}x^4\end{align*}
we see that these cannot be global maximums either.\par
Now to be sure that $0$ is the global maximum value of $s(x,y,z,-x-y-z),$ consider the following
\begin{equation}
s(t(x,y,z,w))=t^4s(x,y,z,w).
\end{equation}
Therefore to prove that $0$ is the global maximum value it is enough to prove that $0$ is the maximum on a neighborhood of the origin. We will consider the closed unit ball in $\mathbb{R}^4$ intersected with our constraint hyperplane: $x+y+z+w=0.$ By homogeneity, if there is point $p\in \{(x,y,z,w):x+y+z+w=0\}$ such that $s(p)>0$ then $s$ will have to achieve its maximum on the boundary of the ball. Assume for the sake of contradiction that this maximum is achieved at $\underline{x}\in \mathbb{S}^4\cap \{(x,y,z,w) : x+y+z+w=0\}.$ 
By homogeneity of the function we get 
\begin{align}
4s(\underline{x})=\sum_i \partial_i s (\underline{x})x_i  
\end{align}
and by Lagrange multipliers we get that there exists $\lambda, \mu\in\mathbb{R}$ such that 
\begin{equation}
\partial_i s(\underline{x})=\lambda+2\mu x_i\quad \forall\hspace{1mm}i\in\{1,2,3,4\}.
\end{equation}
Using $(2.16)$ and $(2.17)$ we can derive
\begin{equation}
4s(\underline{x})=2\mu.
\end{equation}
Suppose that $\underline{x}$ has one coordinate which is $0,$ wlog say $x=0.$ Then
\begin{align*}
\partial_{x} s(x,y,z,w)&=
4x^3-\frac{28}{12}||x||^2x\\
&=0.
\end{align*}
Therefore by $(2.17)$ \begin{equation*}
\lambda=0.
\end{equation*}
This implies that $(2.17)$ can be written as \begin{equation}\partial_is(\underline{x})=2\mu x_i.\end{equation}
Now if $\underline{x}$ has another coordinate equal to zero then the value of $s(\underline{x})$ will be negative by $(2.14),$ so suppose $y,z,w$ are non-zero. Now exactly two of the non-zero coordinates must have the same sign. Assume $y\geq z>0.$ \par
Next consider
 \begin{align*}
\partial_{z} s(x,y,z,w)&=
4z^3-\frac{28}{12}||x||^2z\\
&=4z^3-\frac{28}{12}(x^2+y^2+z^2+w^2)z\\
&\leq 4z^3-\frac{28}{12}z^3-\frac{28}{12}(y^2+z^2+x^2)z\\
&< 0.
\end{align*}
Therefore by $(2.19)$
\begin{equation*}
2\mu z<0,
\end{equation*}
this implies that \begin{equation*}\mu <0\end{equation*}
which is a contradiction to $(2.18).$
We can make a similar argument in the case where two of $y,z,w$ are negative as follows, suppose $0>z\geq w,$ then
 \begin{align*}
\partial_{z} s(x,y,z,w)&=
4z^3-\frac{28}{12}||x||^2z\\
&=4z^3-\frac{28}{12}(x^2+y^2+z^2+w^2)z\\
&\geq 4z^3-\frac{28}{12}z^3-\frac{28}{12}(y^2+z^2+x^2)z\\
&> 0.
\end{align*}
Therefore by $(2.19)$
\begin{equation}
2\mu z>0,
\end{equation}
this implies that 
\begin{equation*}
\mu<0
\end{equation*}
which is a contradiction to $(2.18).$
Therefore $\underline{x}$ can only have non-zero coordinates.
\par

 Suppose that $\underline{x}$ has two positive coordinates and two negative coordinates, wlog assume that the positive coordinates are $x$ and $y$ and assume wlog $x\geq y$ and that the negative coordinates are $z,w$ with $z\geq w.$
Then \begin{align*}
\partial_{y} s(x,y,z,w)&=
4y^3-\frac{28}{12}||x||^2y\\
&=4y^3-\frac{28}{12}(x^2+y^2+z^2+w^2)y\\
&\leq 4y^3-\frac{28}{12}y^3-\frac{28}{12}(y^2+z^2+w^2)y\\
&< 0.
\end{align*}
Therefore \begin{equation*}
\lambda+2\mu y<0,
\end{equation*}
this implies that \begin{equation*}\lambda<0.\end{equation*}
Next consider
 \begin{align*}
\partial_{z} s(x,y,z,w)&=
4z^3-\frac{28}{12}||x||^2z\\
&=4z^3-\frac{28}{12}(x^2+y^2+z^2+w^2)z\\
&\geq 4z^3-\frac{28}{12}z^3-\frac{28}{12}(y^2+z^2+x^2)z\\
&> 0.
\end{align*}
Therefore
\begin{equation*}
\lambda+2\mu z>0,
\end{equation*}
this implies that \begin{equation*}\lambda >0.\end{equation*} Therefore we have reached a contradiction and hence $\underline{x}$ cannot have two positive coordinates and two negative coordinates.
\par
Now we know that $\underline{x}$ must have $3$ coordinates of the same sign and the $4th$ must be of the opposite. 
Suppose $x\geq y\geq z>0$ and $w=-x-y-z.$ Then
\begin{equation}
\partial_x s(\underline{x})=4x^3-\frac{7}{3}x=\lambda+2\mu x
\end{equation}
which implies
\begin{equation}
4x^3-\frac{7}{3}x-2\mu x=\lambda
\end{equation}
the right hand-side of $(2.22)$ is independent of which variable we chose in $(2.21)$ so we get 
\begin{equation}
4x_i^3-\frac{7}{3}x_i-2\mu x_i=4x_j^3-\frac{7}{3}x_j-2\mu x_j.
\end{equation}
It follows that 
\begin{equation}
4(x_i^3-x_j^3)-(\frac{7}{3}+2\mu)(x_i-x_j)=0
\end{equation}
now if $x_i\neq x_j$ then dividing by $x_i-x_j$ gives
\begin{equation}
4(x_i^2+x_ix_j+x_j^2)-(\frac{7}{3}+2\mu)=0
\end{equation}
therefore 
\begin{equation}
4(x_i^2+x_ix_j+x_j^2)=(\frac{7}{3}+2\mu)
\end{equation}
now since the right hand-side is independent of the choice of coordinates we get that if $x, y$ and $z$ are distinct then 
\begin{align}
4(x^2+xy+y^2)&=\frac{7}{3}+2\mu\\
4(x^2+xw+w^2)&=\frac{7}{3}+2\mu
\end{align}
writing $w=-x-y-z$ and taking the average of $(2.27)$ and $(2.28)$ gives
\begin{equation*}
2(x^2-xz+w^2+y^2)=\frac{7}{3}+2\mu.
\end{equation*}
We could have swapped the roles of $x$ and $z$ and gotten 
\begin{equation*}
2(z^2-xz+w^2+y^2)=\frac{7}{3}+2\mu.
\end{equation*}
Which implies 
\begin{equation*}
x=z.
\end{equation*}
Therefore $x,y$ and $z$ cannot be distinct. \par Now suppose two of the coordinates are equal, wlog $x=z,$ then if $y\neq x$ we get the following equations
\begin{align}
4(y^2+yz+z^2)&=\frac{7}{3}+2\mu\\
4(z^2+zw+w^2)&=\frac{7}{3}+2\mu\\
4(y^2+yw+w^2)&=\frac{7}{3}+2\mu
\end{align}
Writing $w=-x-y-z$ and taking the average of equations $(2.29)$ and $(2.30)$ gives
\begin{equation*}2(y^2-zx+w^2+z^2)=\frac{7}{3}+2\mu.
\end{equation*}
Taking the average of $(2.29)$ and $(2.31)$ gives
\begin{equation*}
2(y^2-yx+z^2+w^2)=\frac{7}{3}+2\mu.
\end{equation*}
This implies that $y=z.$\par Therefore $\underline{x}$ is of the form $\alpha \sigma(1,1,1,-3)$ for $\alpha\in\mathbb{R_+}$ and $\sigma\in S_4$ and we get \begin{equation*}
s(\underline{x})=0,\end{equation*}
which is a contradiction.
If $x,y$ and $z$ were negative instead of positive then we would still get a contradiction since the point $-\underline{x}$ would also be a local maximum on which $s$ has the same value by symmetry but by what we've just argued we would get $s(\underline{x})=0.$ Now by corollary $4.4$ of \cite{do2012rotation} we get the result.
\end{proof}

\bibliographystyle{amsplain}
\bibliography{bibliography.bib}

\vfill

\end{document}